\documentclass[12pt,draft]{article}

\usepackage[cp1251]{inputenc}
\usepackage[T2A]{fontenc}
\usepackage[english,ukrainian]{babel}
\usepackage{amsmath,amsfonts,amssymb}
\usepackage{geometry}
\usepackage{cite}
\begin{document}

\begin{flushleft}


\textbf{Tetiana Kasirenko and Aleksandr Murach}\\
\small(Institute of Mathematics, National Academy of Sciences of Ukraine, Kyiv)

\medskip

\large\textbf{ELLIPTIC PROBLEMS WITH BOUNDARY CONDITIONS OF HIGH ORDERS IN H\"ORMANDER SPACES}\normalsize

\medskip\medskip

\textbf{Тетяна Касіренко і Олександр Мурач}\\
\small(Інститут математики НАН України, Київ)

\medskip

\large\textbf{ЕЛІПТИЧНІ ЗАДАЧІ З КРАЙОВИМИ УМОВАМИ ВИСОКИХ ПОРЯДКІВ У ПРОСТОРАХ ХЕРМАНДЕРА}\normalsize

\end{flushleft}

\medskip

\noindent In a class of inner product H\"ormander spaces, we investigate a general elliptic problem for which the maximum of orders of boundary conditions is grater than or equal to the order of elliptic equation. The order of regularity for these spaces is an arbitrary radial positive function RO-varying at infinity in the sense of Avakumovi\'c. We prove that the operator of the problem under investigation is bounded and Fredholm on appropriate pairs of H\"ormander spaces indicated.  A theorem on isomorphism generated by this operator is proved. For generalized solutions to this problem, we establish a local a priory estimate and prove a theorem about their local regularity in H\"ormander spaces. As application, we obtain new sufficient conditions under which given derivatives of the solutions are continuous.

\medskip

\noindent У класі гільбертових просторів Хермандера досліджено загальну еліптичну задачу, для якої максимум порядків крайових умов більший або рівний, ніж порядок еліптичного рівняння. Показником регулярності для цих просторів служить довільна радіальна додатна функція, RO-змінна на нескінченності за Авакумовичем. Доведено, що оператор досліджуваної задачі є обмеженим і нетеровим у підходящих парах вказаних просторів Хермандера. Доведено теорему про ізоморфізм, породжений цим оператором. Для узагальнених розв'язків цієї задачі встановлено локальну апріорну оцінку і доведено теорему про їх локальну регулярність у просторах Хермандера. Як застосування, отримано нові достатні умови неперервності заданих узагальнених похідних розв'язків.

\medskip

\noindent В классе гильбертовых пространств Хермандера исследована общая эллиптическая задача, для которой максимум порядков краевых условий больший или равный, чем порядок эллиптического уравнения. Показателем регулярности для этих пространств служит произвольная радиальная положительная функция, RO-меняющаяся на бесконечности по Авакумовичу. Доказано, что оператор исследуемой задачи является ограниченным и нетеровым в подходящих парах указанных пространств Хермандера. Доказана теорема об  изоморфизме, порожденном этим оператором. Для обобщенных решений этой задачи установлена локальная априорная оценка и доказана теорема об их локальной регулярности в пространствах Хермандера. В качестве приложения получены новые достаточные условия непрерывности заданных обобщенных производных решений.

\bigskip

\noindent\textbf{1. Вступ.} Центральний результат теорії загальних еліптичних крайових задач в обмежених областях з гладкою межею полягає у тому, що ці задачі є нетеровими у підходящих парах функціональних просторів Соболєва або Гельдера (див., наприклад, огляд \cite[\S2]{Agranovich97}, довідник \cite[розд.~III, \S6]{FunctionalAnalysis72} і монографії \cite{Hermander63, LionsMagenes71, Triebel95}). Цей результат має різні застосування, серед яких твердження про підвищення регулярності розв'язків еліптичних задач.  Втім, класичні шкали Соболєва і Гельдера є недостатньо тонко градуйованими для низки задач, що виникають в аналізі і теорії диференціальних рівнянь. У цьому зв'язку Л.~Хермандер \cite{Hermander63} вів широкі класи нормованих функціональних просторів, для яких показником регулярності розподілів служить не число, а досить загальна вагова функція, залежна від частотних змінних. Хермандер \cite{Hermander63, Hermander83} дав застосування цих просторів до дослідження характеру розв'язності і регулярності розв'язків лінійних диференціальних рівнянь з частинними похідними. Простори Хермандера і різні їх узагальнення знайшли  застосування в математичному аналізі, теорії диференціальних рівнянь, теорії випадкових процесів (див. монографії \cite{Jacob010205, MikhailetsMurach14, NicolaRodino10, Paneah00, Triebel01}).

Недавно В.~А.~Михайлець і другий автор цієї статті \cite{MikhailetsMurach05UMJ5, MikhailetsMurach06UMJ2, MikhailetsMurach06UMJ3, MikhailetsMurach07UMJ5, MikhailetsMurach06UMJ11,  MikhailetsMurach08UMJ4, Murach08MFAT2, MikhailetsMurach12BJMA2, MikhailetsMurach14} побудували теорію розв'язності загальних еліптичних систем на гладких многовидах і еліптичних крайових задач у класах гільбертових просторів Хермандера, які отримуються інтерполяцією з функціональним параметром пар гільбертових просторів Соболєва. Показником регулярності для цих просторів служать радіальні функції, правильно змінні на нескінченності за Караматою \cite{Karamata30a} (див. монографії \cite{Seneta76, BinghamGoldieTeugels89}). За допомогою методу інтерполяції з функціональним параметром гільбертових просторів вдалося перенести основні результати “соболєвської” теорії еліптичних рівнянь і задач на зазначені простори Хермандера. Ці результати були доповнені в \cite{Murach09UMJ3, ZinchenkoMurach13UMJ11, ZinchenkoMurach14JMathSci, AnopMurach14MFAT2, AnopMurach14UMJ7, ChepurukhinaMurach15UMJ5, AnopKasirenko16MFAT4} для більш широких класів гільбертових просторів Хермандера. Відмітимо, що згаданий метод інтерполяції виявився плідним і в теорії параболічних початково-крайових задач \cite{LosMikhailetsMurach17CPAA, LosMurach17OpenMath}.

У побудованій теорії розглянуто виключно еліптичні задачі, у яких порядки крайових умов менші, ніж порядок еліптичного рівняння. Мета даної статті~--- доповнити цю теорію результатами про характер розв'язності і властивості розв'язків еліптичних задач, у яких порядок принаймні однієї з крайових умов більший або рівний за порядок еліптичного рівняння. Ці задачі будемо досліджувати у класі гільбертових просторів Хермандера, показником регулярності для яких служить довільна радіальна функція, RO-змінна на нескінченності за Авакумовичем \cite{Avakumovic36} (див. монографію \cite{Seneta76}). Цей клас був виділений в \cite{MikhailetsMurach09Dop3, MikhailetsMurach13UMJ3} і названий розширеною соболєвською шкалою. Він містить уточнену соболєвську шкалу та складається з усіх гільбертових просторів, інтерполяційних відносно пар гільбертових просторів Соболєва.

Робота складається з 7 пунктів. Пункт~1 є вступом. У п.~2 сформульовано еліптичну крайову задачу, яка досліджується, і розглянуто формально спряжену до неї задачу відносно спеціальної формули Гріна. У п.~3 наведено означення функціональних просторів Хермандера, які утворюють розширену соболєвську шкалу. Пункт~4 містить основні результати роботи про властивості досліджуваної задачі в просторах Хермандера. У п.~5, як застосування основних результатів, отримано достатні умови неперервності узагальнених похідних розв'язків досліджуваної задачі, зокрема, умови класичності її узагальненого розв'язку. Пункт~6 присвячений інтерполяції з функціональним параметром пар гільбертових просторів та її застосуванням до розширеної соболєвської шкали. Результати роботи доведено у заключному п.~7.

\textbf{2. Постановка задачі.} Нехай $\Omega$~--- довільна обмежена область у евклідовому просторі $\mathbb{R}^{n}$, де $n\geq2$. Припускаємо, що межа $\Gamma$ цієї області є нескінченно гладким компактним многовидом вимірності $n-1$, причому $C^\infty$-структура на  $\Gamma$ індукована простором~$\mathbb{R}^{n}$.

В області $\Omega$ розглянемо таку крайову задачу:
\begin{gather}\label{1f1}
Au=f\quad\mbox{в}\quad\Omega,\\
B_{j}u=g_{j}\quad\mbox{на}\quad\Gamma,
\quad j=1,...,q. \label{1f2}
\end{gather}
Тут
$$
A:=A(x,D):=\sum_{|\mu|\leq 2q}a_{\mu}(x)D^{\mu}
$$
є лінійним диференціальним оператором на $\overline{\Omega}:=\Omega\cup\Gamma$ довільного парного порядку $2q\geq\nobreak2$, а кожне
$$
B_{j}:=B_{j}(x,D):=\sum_{|\mu|\leq m_{j}}b_{j,\mu}(x)D^{\mu}
$$
є крайовим лінійним диференціальним оператором на $\Gamma$ довільного порядку $m_{j}\geq0$. Усі коефіцієнти $a_{\mu}(x)$ і $b_{j,\mu}(x)$ цих диференціальних операторів є нескінченно гладкими комплекснозначними функціями на $\overline{\Omega}$ і $\Gamma$ відповідно. Взагалі, у роботі функції та розподіли припускаються комплекснозначними і тому усі розглянуті функціональні простори вважаються комплексними.

У наведених формулах і далі використано такі стандартні позначення:
$\mu:=(\mu_{1},\ldots,\mu_{n})$~--- мультиіндекс з невід'ємними цілими компонентами, $|\mu|:=\mu_{1}+\cdots+\mu_{n}$,
$D^{\mu}:=D_{1}^{\mu_{1}}\ldots D_{n}^{\mu_{n}}$, де $D_{k}:=i\partial/\partial x_{k}$ для кожного номера $k\in\{1,...,n\}$, $i$~--- уявна одиниця, а $x=(x_1,\ldots,x_n)$~--- довільна точка простору $\mathbb{R}^{n}$. Також покладемо $D_{\nu}:=i\partial/\partial\nu$, де $\nu(x)$~--- орт внутрішньої нормалі до межі $\Gamma$ у точці $x\in\Gamma$.

У роботі припускаємо, що крайова задача \eqref{1f1}, \eqref{1f2}
є еліптичною в області $\Omega$, тобто диференціальний оператор $A$ є правильно еліптичним на $\overline{\Omega}$, а набір $B:=(B_{1},\ldots,B_q)$  крайових диференціальних операторів задовольняє умову Лопатинського щодо $A$ на $\Gamma$ (див., наприклад, огляд \cite[п.~1.2]{Agranovich97} або довідник \cite[розд.~III, \S~6, пп. 1, 2]{FunctionalAnalysis72}).

\textbf{Приклад 1.} Розглянемо крайову задачу, яка складається з диференціального рівняння \eqref{1f1}, де диференціальний оператор $A$ правильно еліптичний на $\overline{\Omega}$, і крайових умов
\begin{equation*}
\frac{\partial^{k+j-1}u}{\partial\zeta^{k+j-1}}+\sum_{|\mu|< k+j-1}b_{j,\mu}(x)D^{\mu}=g_j\quad\text{на}\quad\Gamma,
\quad j=1,...,q.
\end{equation*}
Тут ціле число $k\geq0$, а $\zeta:\Gamma\to\mathbb{R}^{n}$ є нескінченно гладким полем векторів $\zeta(x)$, недотичних до $\Gamma$ у точці $x\in\Gamma$. Безпосередньо перевіряється, що ця крайова задача є еліптичною в області $\Omega$. Якщо $0\leq k\leq q$, то вона є регулярною еліптичною  (див., наприклад, \cite[п.~5.2.1, зауваження~4]{Triebel95}). Важливий окремий випадок цієї задачі отримуємо, поклавши $A:=\Delta^{q}$, де $\Delta$~--- оператор Лапласа, та $\zeta(x):=\nu(x)$ для усіх $x\in\Gamma$.

Надалі припускаємо, що
\begin{equation*}
m:=\max\{m_{1},\ldots,m_{q}\}\geq2q.
\end{equation*}

Пов'яжемо із задачею \eqref{1f1}, \eqref{1f2} лінійне відображення
\begin{equation}\label{1f3}
u\mapsto(Au,Bu)=(Au,B_{1}u,\ldots,B_{q}u),\quad\mbox{де}\quad u\in C^{\infty}(\overline{\Omega}).
\end{equation}
Мета роботи~--- дослідити властивості продовження за неперервністю цього відображення у підходящих парах функціональних просторів Хермандера.

Для опису області значень цього продовження нам потрібна така спеціальна формула Гріна \cite[формула (4.1.10)]{KozlovMazyaRossmann97}:
\begin{gather*}
(Au,v)_{\Omega}+\sum_{j=1}^{m-2q+1}(D_{\nu}^{j-1}Au,w_{j})_{\Gamma}+
\sum_{j=1}^{q}(B_{j}u,h_{j})_{\Gamma}=\\
=(u,A^{+}v)_{\Omega}+\sum_{k=1}^{m+1}\biggl(D_{\nu}^{k-1}u,K_{k}v+
\sum_{j=1}^{m-2q+1}R_{j,k}^{+}w_{j}+
\sum_{j=1}^{q}Q_{j,k}^{+}h_{j}\biggr)_{\Gamma},
\end{gather*}
де $u,v\in C^{\infty}(\overline{\Omega})$, $w_{1},\ldots,w_{m-2q+1},h_{1},\ldots,h_{q}
\in C^{\infty}(\Gamma)$ та через $(\cdot,\cdot)_{\Omega}$ і $(\cdot,\cdot)_{\Gamma}$ позначено відповідно скалярні добутки у гільбертових просторах $L_{2}(\Omega)$ і $L_{2}(\Gamma)$ функцій квадратично інтегровних на $\Omega$ і $\Gamma$ відносно мір Лебега. Тут $A^{+}$~--- диференціальний оператор, формально спряжений до $A$, тобто
$$
(A^{+}v)(x):=\sum_{|\mu|\leq2q}D^{\mu}(\overline{a_{\mu}(x)}v(x)).
$$
Окрім того, усі $R_{j,k}^{+}$ і $Q_{j,k}^{+}$ є дотичними диференціальними операторами, формально спряженими  відповідно до $R_{j,k}$ і $Q_{j,k}$ відносно $(\cdot,\cdot)_{\Gamma}$, а дотичні лінійні диференціальні оператори $R_{j,k}:=R_{j,k}(x,D_{\tau})$ і $Q_{j,k}:=Q_{j,k}(x,D_{\tau})$ узяті із зображення крайових диференціальних операторів $D_{\nu}^{j-1}A$ і $B_{j}$ у вигляді
\begin{gather*}
D_{\nu}^{j-1}A(x,D)=\sum_{k=1}^{m+1}R_{j,k}(x,D_{\tau})D_{\nu}^{k-1},\quad
j=1,\ldots,m-2q+1,\\
B_{j}(x,D)=\sum_{k=1}^{m+1}Q_{j,k}(x,D_{\tau})D_{\nu}^{k-1},\quad
j=1,\ldots,q.
\end{gather*}
Відмітимо, що $\mathrm{ord}\,R_{j,k}\leq 2q+j-k$ і $\mathrm{ord}\,Q_{j,k}\leq m_{j}-k+1$, причому, звісно, $R_{j,k}=0$ при $k\geq2q+j+1$ і $Q_{j,k}=0$ при $k\geq m_{j}+2$. Нарешті, кожне $K_{k}:=K_{k}(x,D)$~--- деякий крайовий лінійний диференціальний оператор на $\Gamma$ порядку $\mathrm{ord}\,K_{k}\leq2q-k$ з коефіцієнтами класу $C^{\infty}(\overline{\Omega})$.

Спеціальна формула Гріна приводить до такої крайової
задачі в області $\Omega$:
\begin{gather}\label{1f4}
A^{+}v=\omega\quad\mbox{в}\quad\Omega,\\
K_{k}v+\sum_{j=1}^{m-2q+1}R_{j,k}^{+}w_{j}+
\sum_{j=1}^{q}Q_{j,k}^{+}h_{j}=\theta_{k}\quad
\mbox{на}\quad\Gamma,\quad k=1,...,m+1. \label{1f5}
\end{gather}
Ця задача містить, окрім невідомої функції $v$ на $\Omega$, ще
$m-q+1$ додаткових невідомих функцій $w_{1},\ldots,w_{m-2q+1},h_{1},\ldots,h_{q}$ на межі $\Gamma$. Задачу \eqref{1f4}, \eqref{1f5} називають формально спряженою до задачі \eqref{1f1}, \eqref{1f2} відносно розглянутої спеціальної формули Гріна. Відомо \cite[теорема~4.1.1]{KozlovMazyaRossmann97}, що крайова задача \eqref{1f1}, \eqref{1f2} еліптична тоді і тільки тоді, коли формально спряжена задача \eqref{1f4}, \eqref{1f5} еліптична як крайова задача з додатковими невідомими функціями на межі області.

\textbf{Приклад 2.} Запишемо спеціальну формулу Гріна для еліптичної крайової задачі
\begin{equation}\label{1f1ex2}
\Delta u=f\;\;\text{в}\;\;\Omega,\qquad
\frac{\partial^{2}u}{\partial\nu^{2}}=g\;\;\text{на}\;\;\Gamma,
\end{equation}
заданої в крузі $\Omega:=\{(x_1,x_2)\in \mathbb{R}^{2}: x_{1}^2+x_{2}^2<1\}$. Відмітимо, що $\Delta u=\partial_{\nu}^{2}u-\partial_{\nu}u+\partial_{\varphi}^{2}u$ на~$\Gamma$; тут $\partial_{\nu}:=\partial/\partial\nu=-\partial/\partial\varrho$ і $\partial_{\varphi}:=\partial/\partial\varphi$, а $(\varrho,\varphi)$~--- полярні координати. Застосувавши другу класичну формулу Гріна для оператора Лапласа, отримаємо, що
\begin{gather*}
(\Delta u,v)_{\Omega}+(\Delta u,w)_{\Gamma}+
(\partial_{\nu}^{2}u,h)_{\Gamma}=\\
=(u,\Delta v)_{\Omega}-(\partial_{\nu}u,v)_{\Gamma}+
(u,\partial_{\nu}v)_{\Gamma}+(\partial_{\nu}^{2}u-
\partial_{\nu}u+\partial_{\varphi}^{2}u,w)_{\Gamma}+
(\partial_{\nu}^{2}u,h)_{\Gamma}=\\
=(u,\Delta v)_{\Omega}+
(u,\partial_{\nu}v+\partial_{\varphi}^{2}w)_{\Gamma}+
(\partial_{\nu}u,-v-w)_{\Gamma}+
(\partial_{\nu}^{2}u,w+h)_{\Gamma}
\end{gather*}
для довільних функцій $u,v\in C^{\infty}(\overline{\Omega})$ і $w,h\in C^{\infty}(\Gamma)$. Отже, спеціальна формула Гріна для крайової задачі \eqref{1f1ex2} набирає вигляду
\begin{gather*}
(\Delta u,v)_{\Omega}+(\Delta u,w)_{\Gamma}+
(\partial_{\nu}^{2}u,h)_{\Gamma}=\\
=(u,\Delta v)_{\Omega}+
(u,\partial_{\nu}v+\partial_{\varphi}^{2}w)_{\Gamma}+
(D_{\nu}u,-iv-iw)_{\Gamma}+(D_{\nu}^{2}u,-w-h)_{\Gamma}.
\end{gather*}
Тому крайова задача
\begin{gather*}
\Delta v=\omega\quad\mbox{в}\quad\Omega,\\
\partial_{\nu}v+\partial_{\varphi}^{2}w=\theta_{1},\quad
-iv-iw=\theta_{2},\quad
-w-h=\theta_{3}\quad\mbox{на}\quad\Gamma
\end{gather*}
є формально спряженою до задачі \eqref{1f1ex2} відносно цієї формули Гріна. Отримана формально спряжена задача містить дві додаткові невідомі функції $w$ і $h$ на $\Gamma$.

\textbf{3. Простори Хермандера і розширена соболєвська шкала.} Еліптичну крайову задачу \eqref{1f1}, \eqref{1f2} будемо досліджувати у підходящих парах гільбертових просторів Хермандера \cite[п.~2.2]{Hermander63}, які утворюють розширену соболєвську шкалу, введену в \cite{MikhailetsMurach09Dop3, MikhailetsMurach13UMJ3}. Нагадаємо означення цих просторів і деякі їх властивості, потрібні у подальшому.

Для просторів Хермандера, які використовуються у роботі, показником регулярності розподілів служить функціональний параметр $\alpha\in\mathrm{RO}$. За означенням, клас $\mathrm{RO}$ складається з усіх вимірних за Борелем функцій
$\alpha:\nobreak[1,\infty)\rightarrow(0,\infty)$, для яких існують числа $b>1$ і $c\geq1$ такі, що $c^{-1}\leq\alpha(\lambda t)/\alpha(t)\leq c$ для довільних $t\geq1$ і $\lambda\in[1,b]$
(сталі $b$ і $c$ можуть залежати від $\alpha$). Такі
функції називають RO-змінними на нескінченності. Клас RO введений
В.~Г.~Авакумовичем \cite{Avakumovic36} у 1936~р. і достатньо повно вивчений (див., наприклад, монографії \cite[додаток~1]{Seneta76} і \cite[пп.~2.0~-- 2.2]{BinghamGoldieTeugels89}).

Цей клас допускає простий опис, а саме:
$$
\alpha\in\mathrm{RO}\;\;\Leftrightarrow\;\;\alpha(t)=\exp\Biggl(\beta(t)+
\int\limits_{1}^{\:t}\frac{\gamma(\tau)}{\tau}\;d\tau\Biggr)\;\,
\mbox{для}\;\,t\geq1,
$$
де дійсні функції $\beta$ і $\gamma$ вимірні за Борелем і обмежені на півосі $[1,\infty)$ (див., наприклад, \cite[додаток~1, теорема~1]{Seneta76}).

Для нас важлива така властивість класу $\mathrm{RO}$: для
кожної функції $\alpha\in\mathrm{RO}$ існують числа $s_{0},s_{1}\in\mathbb{R}$,
$s_{0}\leq s_{1}$, і $c_{0},c_{1}>0$ такі, що
\begin{equation}\label{Ax=b3}
c_{0}\lambda^{s_{0}}\leq\frac{\alpha(\lambda t)}{\alpha(t)}\leq
c_{1}\lambda^{s_{1}} \quad\mbox{для всіх}\quad t\geq1,\;\;\lambda\geq1
\end{equation}
(див. \cite[додаток~1, теорема~2]{Seneta76}). Покладемо
\begin{gather*}
\sigma_{0}(\alpha):=
\sup\,\{s_{0}\in\mathbb{R}:\,\mbox{виконується ліва нерівність в \eqref{Ax=b3}}\},\\
\sigma_{1}(\alpha):=\inf\,\{s_{1}\in\mathbb{R}:\,\mbox{виконується права нерівність
в \eqref{Ax=b3}}\}.
\end{gather*}
Числа
$\sigma_{0}(\alpha)$ і $\sigma_{1}(\alpha)$ є відповідно нижнім і верхнім індексами Матушевської \cite{Matuszewska64} функції $\alpha\in\mathrm{RO}$ (див. також монографію \cite[п.~2.1.2]{BinghamGoldieTeugels89}). Звісно, $-\infty<\sigma_{0}(\alpha)\leq\sigma_{1}(\alpha)<\infty$.

Наведемо деякі характерні приклади функцій, RO-змінних на нескінченності.

\textbf{Приклад 3.} Розглянемо неперервну функцію $\alpha:[1,\infty)\rightarrow(0,\infty)$ таку, що
\begin{equation*}
\alpha(t):=t^{s}(\ln t)^{r_{1}}(\ln\ln
t)^{r_{2}}\ldots(\underbrace{\ln\ldots\ln}_{k\;\mbox{\small
разів}} t)^{r_{k}}\quad\mbox{при}\quad t\gg1.
\end{equation*}
Тут довільно вибрано ціле число $k\geq1$ і дійсні числа $s,r_{1},\ldots,r_{k}$. Функція $\alpha$ належить до класу $\mathrm{RO}$ і для неї $\sigma_{0}(\alpha)=\sigma_{1}(\alpha)=s$.

Взагалі, до класу $\mathrm{RO}$ належить будь-яка вимірна функція
$\alpha:[1,\infty)\rightarrow(0,\infty)$, яка обмежена і відокремлена від нуля на кожному компакті і є правильно змінною на нескінченності за Й.~Караматою \cite{Karamata30a}. Остання властивість значить, що $\alpha(\lambda t)/\alpha(t)\to\lambda^{s}$ при $t\to\infty$ для деякого $s\in\mathbb{R}$. Індекси Матушевської такої функції дорівнюють числу $s$, яке називають порядком змінення функції на нескінченності. Правильно змінні функції широко застосовуються у математиці (див. монографії \cite{Seneta76, BinghamGoldieTeugels89}).

\textbf{Приклад 4.} Нехай $\theta\in\mathbb{R}$, $\delta>0$ і $r\in(0,1]$. Покладемо
\begin{equation*}
\alpha(t):=\left\{
\begin{array}{ll}
t^{\theta+\delta\sin(\ln\ln t)^{r}}\; &\hbox{при}\;t>e,\\
t^{\theta}\; &\hbox{при}\;1\leq t\leq e.
\end{array}\right.
\end{equation*}
Тоді $\alpha\in\mathrm{RO}$, причому $\sigma_{0}(\alpha)=\theta-\delta$ і $\sigma_{1}(\alpha)=\theta+\delta$ \cite[приклад~6]{Chepuruhina15Coll2}.

Нехай $\alpha\in\mathrm{RO}$. Дамо означення простору Хермандера $H^{\alpha}$ спочатку на $\mathbb{R}^{n}$, де ціле $n\geq1$, а потім на $\Omega$ і $\Gamma$. Цей простір складається з розподілів (узагальнених функцій), які нам зручно трактувати як \emph{анти}лінійні функціонали на відповідному просторі основних функцій.

За означенням, лінійний простір $H^{\alpha}(\mathbb{R}^{n})$ складається з усіх повільно зростаючих на $\mathbb{R}^{n}$ розподілів $w$ таких, що їх перетворення Фур'є $\widehat{w}$ локально інтегровне за Лебегом на $\mathbb{R}^{n}$ і задовольняє умові
$$
\int\limits_{\mathbb{R}^{n}}
\alpha^2(\langle\xi\rangle)\,|\widehat{w}(\xi)|^2\,d\xi
<\infty,
$$
де $\langle\xi\rangle:=(1+|\xi|^{2})^{1/2}$ є згладженим модулем вектора $\xi\in\mathbb{R}^{n}$. Цей простір наділений скалярним добутком
$$
(w_{1},w_{2})_{H^{\alpha}(\mathbb{R}^{n})}:=
\int\limits_{\mathbb{R}^{n}}
\alpha^2(\langle\xi\rangle)\,
\widehat{w_{1}}(\xi)\,\overline{\widehat{w_{2}}(\xi)}\,d\xi
$$
і відповідною нормою
$$
\|w\|_{H^{\alpha}(\mathbb{R}^{n})}:=
(w,w)_{H^{\alpha}(\mathbb{R}^{n})}^{1/2}
$$
та є гільбертовим і сепарабельним відносно цієї норми.

Простір $H^{\alpha}(\mathbb{R}^{n})$~--- гільбертів ізотропний
випадок просторів $\mathcal{B}_{p,k}$, введених і досліджених
Л.~Хермандером в \cite[п. 2.2]{Hermander63} (див також його
монографію \cite[п. 10.1]{Hermander83}). А саме,
$H^{\alpha}(\mathbb{R}^{n})=\mathcal{B}_{p,k}$, якщо $p=2$ і
$k(\xi)=\alpha(\langle\xi\rangle)$ при $\xi\in\mathbb{R}^{n}$.
Зауважимо, що у гільбертовому випадку $p=2$ простори Хермандера
збігаються з просторами, введеними Л.~Р.~Волевичем і
Б.~П.~Панеяхом \cite[\S~2]{VolevichPaneah65}.

Якщо $\alpha(t)\equiv t^{s}$ для деякого $s\in\mathbb{R}$, то
$H^{\alpha}(\mathbb{R}^{n})=:H^{(s)}(\mathbb{R}^{n})$ є
гільбертів простір Соболєва порядку $s$. Взагалі,
\begin{equation}\label{1f7}
s_{0}<\sigma_{0}(\alpha)\leq\sigma_{1}(\alpha)<s_{1}\;\Rightarrow\;
H^{(s_1)}(\mathbb{R}^{n})\hookrightarrow
H^{\alpha}(\mathbb{R}^{n})\hookrightarrow
H^{(s_0)}(\mathbb{R}^{n}),
\end{equation}
причому обидва вкладення неперервні й щільні.

Cлідуючи \cite{MikhailetsMurach13UMJ3, MikhailetsMurach14}, клас функціональних просторів
$\{H^{\alpha}(\mathbb{R}^{n}):\alpha\in\mathrm{RO}\}$
називаємо розширеною соболєвською шкалою на $\mathbb{R}^{n}$. Її аналоги для $\Omega$ і $\Gamma$ будуються стандартним чином (див. \cite[с.~4]{MikhailetsMurach15ResMath1} і \cite[с.~30]{MikhailetsMurach09Dop3}). Наведемо відповідні означення; тепер $n\geq2$.

За означенням, лінійний простір $H^{\alpha}(\Omega)$ складається зі звужень в область $\Omega$ всіх розподілів $w\in H^{\alpha}(\mathbb{R}^{n})$ і наділений нормою
$$
\|v\|_{H^{\alpha}(\Omega)}:=
\inf\bigl\{\,\|w\|_{H^{\alpha}(\mathbb{R}^{n})}:\,
w\in H^{\alpha}(\mathbb{R}^{n}),\ w=v\;\,\mbox{в}\;\,\Omega\,\bigr\},
$$
де $v\in H^{\alpha}(\Omega)$. Простір $H^{\alpha}(\Omega)$
гільбертів і сепарабельний відносно цієї норми, а множина $C^{\infty}(\overline{\Omega})$ щільна в ньому.

Лінійний простір $H^{\alpha}(\Gamma)$ складається, коротко кажучи, з усіх розподілів на $\Gamma$, які в локальних координатах дають елементи простору $H^{\alpha}(\mathbb{R}^{n-1})$. Дамо детальне означення. Довільно виберемо скінченний атлас із $C^{\infty}$-структури на многовиді $\Gamma$, утворений локальними картами $\pi_j:\mathbb{R}^{n-1}\leftrightarrow \Gamma_{j}$, де
$j=1,\ldots,\varkappa$. Тут відкриті множини $\{\Gamma_{1},\ldots,\Gamma_{\varkappa}\}$ складають покриття
многовиду $\Gamma$. Виберемо також функції $\chi_j\in C^{\infty}(\Gamma)$, де $j=1,\ldots,\varkappa$, які утворюють розбиття одиниці на $\Gamma$, що задовольняє умову $\mathrm{supp}\,\chi_j\subset \Gamma_j$.

За означенням, лінійний простір $H^{\alpha}(\Gamma)$ складається з усіх розподілів $h$ на $\Gamma$ таких, що $(\chi_{j}h)\circ\pi_{j}\in H^{\alpha}(\mathbb{R}^{n-1})$ для усіх $j\in\{1,\ldots,\varkappa\}$. Тут $(\chi_{j}h)\circ\pi_{j}$ є представленням
розподілу $h$ у локальній карті $\pi_{j}$. Простір $H^{\alpha}(\Gamma)$ наділений нормою
$$
\|h\|_{H^{\alpha}(\Gamma)}:=\biggl(\sum_{j=1}^{\varkappa}\,
\|(\chi_{j}h)\circ\pi_{j}\|_
{H^{\alpha}(\mathbb{R}^{n-1})}^{2}\biggr)^{1/2}.
$$
Він гільбертів і сепарабельний відносно цієї норми та з точністю до еквівалентності норм не залежить від зробленого вибору атласу і розбиття одиниці \cite[с.~32]{MikhailetsMurach09Dop3}. Множина $C^{\infty}(\Gamma)$ щільна в $H^{\alpha}(\Gamma)$.

Щойно означені функціональні простори утворюють розширені соболєвські шкали $\{H^{\alpha}(\Omega):\alpha\in\mathrm{RO}\}$ і
$\{H^{\alpha}(\Gamma):\alpha\in\mathrm{RO}\}$ на $\Omega$ і $\Gamma$ відповідно. Вони містять гільбертові шкали просторів Соболєва: якщо $\alpha(t)\equiv t^{s}$ для деякого $s\in\mathbb{R}$, то $H^{\alpha}(\Omega)=:H^{(s)}(\Omega)$ і
$H^{\alpha}(\Gamma)=:H^{(s)}(\Gamma)$ є гільбертовими просторами Соболєва порядку $s$.

Відмітимо таку властивість цих шкал, яка випливає з \cite[теореми 2.2.2, 2.2.3]{Hermander63}. Нехай $\alpha,\eta\in\mathrm{RO}$ і $\Lambda\in\{\Omega,\Gamma\}$. Функція $\alpha/\eta$ обмежена в околі нескінченності тоді і тільки тоді, коли $H^{\eta}(\Lambda)\hookrightarrow H^\alpha(\Lambda)$. Це вкладення неперервне і щільне. Воно компактне тоді і тільки тоді, коли
$\alpha(t)/\eta(t)\rightarrow0$ при $t\rightarrow\infty$. Зокрема, виконуються властивість \eqref{1f7}, якщо у ній замінити $\mathbb{R}^{n}$ на $\Omega$ або $\Gamma$, при цьому вкладення будуть компактними і щільними.

\textbf{4. Основні результати.} Сформулюємо наші результати про властивості еліптичної крайової задачі \eqref{1f1}, \eqref{1f2} у просторах Хермандера $H^{\alpha}$, розглянутих вище. Для них показник регулярності матиме вигляд $\alpha(t)\equiv\varphi(t)t^{s}$, де $\varphi\in\mathrm{RO}$ і $s\in\mathbb{R}$. Для того, щоб не вказувати аргумент $t$ у показнику будемо використовувати функціональний параметр $\varrho(t):=t$ аргументу $t\geq1$ й записувати $\alpha$ у вигляді $\varphi\varrho^s$. Якщо $\varphi\in\mathrm{RO}$, то, звісно, $\varphi\varrho^s\in\mathrm{RO}$ та $\sigma_j(\varphi\varrho^s)=\sigma_j(\varphi)+s$ для кожного $j\in\{0,1\}$.

Позначимо через $N$ лінійний простір усіх розв'язків
$u\in C^{\infty}(\overline{\Omega})$ крайової задачі \eqref{1f1}, \eqref{1f2} у випадку, коли $f=0$ в $\Omega$ і кожне $g_{j}=0$ на~$\Gamma$. Позначимо також через $N_{\star}$ лінійний простір усіх розв'язків
$$
(v,w_{1},\ldots,w_{m-2q+1},h_{1},\ldots,h_{q})\in C^{\infty}(\overline{\Omega})\times(C^{\infty}(\Gamma))^{m-q+1}
$$
формально спряженої крайової задачі \eqref{1f4}, \eqref{1f5} у випадку, коли $\omega=0$ в $\Omega$ і кожне $\theta_{k}=0$ на~$\Gamma$. Оскільки обидві задачі еліптичні в $\Omega$, то простори $N$ і $N_{\star}$ скінченновимірні \cite[наслідок~4.1.1]{KozlovMazyaRossmann97}.

\textbf{Теорема 1.} \it Нехай $\varphi\in\mathrm{RO}$ і
$\sigma_0(\varphi)>m+1/2$. Тоді відображення \eqref{1f3} продовжується єдиним чином (за неперервністю) до обмеженого оператора
\begin{equation}\label{1f8}
(A,B):H^{\varphi}(\Omega)\rightarrow
H^{\varphi\varrho^{-2q}}(\Omega)\oplus\bigoplus_{j=1}^{q}
H^{\varphi\varrho^{-m_j-1/2}}(\Gamma)=:
\mathcal{H}^{\varphi\varrho^{-2q}}(\Omega,\Gamma).
\end{equation}
Цей оператор нетерів. Його ядро дорівнює $N$, а область значень складається з усіх векторів $(f,g_1,\ldots,g_q)\in\mathcal{H}^{\varphi\varrho^{-2q}}(\Omega,\Gamma)$ таких, що
\begin{equation}\label{1f9}
\begin{gathered}
(f,v)_\Omega+\sum_{j=1}^{m-2q+1}(D_{\nu}^{j-1}f,w_{j})_{\Gamma}+
\sum_{j=1}^{q}(g_j,h_{j})_{\Gamma}=0 \\
\mbox{для всіх} \quad (v,w_{1},\ldots,w_{m-2q+1},h_{1},\ldots,h_{q})\in N_{\star}.
\end{gathered}
\end{equation}
Індекс оператора \eqref{1f8} дорівнює $\dim N-\dim N_{\star}$ та не залежить від $\varphi$.\rm

Як і раніше, у формулі \eqref{1f9} через $(\cdot,\cdot)_{\Omega}$ і $(\cdot,\cdot)_{\Gamma}$ позначено скалярні добутки у гільбертових просторах $L_{2}(\Omega)$ і $L_{2}(\Gamma)$ відповідно. Тут згідно з твердженням~4, поданим нижче у п.~6, для кожної функції $f\in H^{\varphi\varrho^{-2q}}(\Omega)$, де $\sigma_0(\varphi)>m+1/2$, коректно означені образи
$$
D_{\nu}^{j-1}f\in H^{\varphi\varrho^{-2q-j+1/2}}(\Gamma)\subset L_{2}(\Gamma)
$$
відносно крайового оператора $D_{\nu}^{j-1}$ порядку $j-1\leq m-2q$.

У зв'язку з теоремою~1 нагадаємо, що лінійний обмежений оператор $T:E_{1}\rightarrow E_{2}$, де $E_{1}$ і $E_{2}$~--- банахові
простори, називають нетеровим, якщо його ядро $\ker T$ і коядро $E_{2}/T(E_{1})$ скінченновимірні. Якщо цей оператор нетерів, то його область значень замкнена в просторі $E_{2}$ (див., наприклад, \cite[Лемма~19.1.1]{Hermander85}) і для нього означений скінченний індекс
$$
\mathrm{ind}\,T:=\dim\ker T-\dim(E_{2}/T(E_{1})).
$$

Зокрема, для еліптичної крайової задачі з прикладу~2 безпосередньо перевіряється, що $\dim N=\dim N_{\star}=3$ і тому індекс оператора \eqref{1f8} дорівнює нулю.

Відмітимо, що умову $\sigma_0(\varphi)>m+1/2$ у теоремі~1 не можна відкинути чи послабити. Зокрема, якщо $\varphi(t)\equiv t^{s}$ для деяких дійсного $s\leq m_{j}+1/2$ і цілого $j\in\{1,\ldots,q\}$, то відображення $u\mapsto B_{j}u$, де $u\in C^{\infty}(\overline{\Omega})$, не можна продовжити до неперервного лінійного оператора, що діє з простору Соболєва $H^{(s)}(\Omega)$ у лінійний топологічний простір $\mathcal{D}'(\Gamma)$ усіх розподілів на~$\Gamma$ (див., наприклад, \cite[Зауваження~3.5]{MikhailetsMurach14}).

У випадку, коли $N=\{0\}$ і $N_{\star}=\{0\}$, оператор \eqref{1f8} здійснює ізоморфізм між просторами $H^{\varphi}(\Omega)$ і $\mathcal{H}^{\varphi\varrho^{-2q}}(\Omega,\Gamma)$. Це випливає з теореми 1 і теореми Банаха про обернений оператор. У загальній ситуації оператор \eqref{1f8} породжує ізоморфізм між деякими їх підпросторами скінченної ковимірності. Ці підпростори і проектори на них зручно будувати у такий спосіб.

Розглянемо розклад простору $H^{\varphi}(\Omega)$, де $\sigma_0(\varphi)>0$, у пряму суму підпросторів
\begin{gather}\label{1f10}
H^{\varphi}(\Omega)=N\dotplus\bigl\{u\in
H^{\varphi}(\Omega):\,(u,w)_\Omega=0\;\;\mbox{для всіх}
\;\;w\in N\bigr\}.
\end{gather}
Ця рівність правильна, оскільки вона є звуженням розкладу простору $L_{2}(\Omega)$ в ортогональну суму підпростору $N$ і його доповнення.
Стосовно розкладу простору $\mathcal{H}^{\varphi\varrho^{-2q}}(\Omega,\Gamma)$ скористаємося таким результатом.

\textbf{Лема 1.} \it Існує скінченновимірний простір $G\subset C^{\infty}(\overline{\Omega})\times(C^{\infty}(\Gamma))^{q}$ такий, що для кожного $\varphi\in\mathrm{RO}$ з $\sigma_0(\varphi)>m+1/2$ правильний розклад
простору $\mathcal{H}^{\varphi\varrho^{-2q}}(\Omega,\Gamma)$ у пряму суму підпросторів
\begin{gather}\label{1f11}
\mathcal{H}^{\varphi\varrho^{-2q}}(\Omega,\Gamma)=G\dotplus
\bigl\{(f,g_1,\ldots,g_q)\in
\mathcal{H}^{\varphi\varrho^{-2q}}(\Omega,\Gamma):
\mbox{виконується \eqref{1f9}}\bigr\};
\end{gather}
при цьому $\dim G=\dim N_{\star}$. \rm

Позначимо через $P$ і  $Q$ косі проектори відповідно просторів $H^{\varphi}(\Omega)$ і $\mathcal{H}^{\varphi\varrho^{-2q}}(\Omega,\Gamma)$ на другі
доданки в сумах \eqref{1f10} і \eqref{1f11} паралельно першим доданкам. Звісно, ці проектори не залежать від $\varphi$.

\textbf{Теорема 2.} \it Нехай $\varphi\in\mathrm{RO}$ і
$\sigma_0(\varphi)>m+1/2$. Тоді звуження відображення \eqref{1f8} на підпростір $P(H^{\varphi}(\Omega))$ є ізоморфізмом
\begin{equation}\label{1f12}
(A,B):\,P(H^{\varphi}(\Omega))\leftrightarrow
Q(\mathcal{H}^{\varphi\varrho^{-2q}}(\Omega,\Gamma)).
\end{equation} \rm

Дослідимо властивості узагальнених розв'язків еліптичної крайової задачі \eqref{1f1}, \eqref{1f2} у просторах Хермандера. Нагадаємо означення таких розв'язків. Покладемо
$$
H^{m+1/2+}(\Omega):=
\bigcup_{\substack{\alpha\in\mathrm{RO}:\\\sigma_{0}(\alpha)>m+1/2}}
H^{\alpha}(\Omega)=\bigcup_{s>m+1/2}H^{(s)}(\Omega);
$$
тут остання рівність правильна з огляду на властивість \eqref{1f7}. Згідно з теоремою 1, для кожної функції $u\in H^{m+1/2+}(\Omega)$ коректно означений вектор
\begin{equation*}
(f,g):=(f,g_{1},\ldots,g_{q}):=(A,B)u
\in L_{2}(\Omega)\times(L_{2}(\Gamma))^{q}.
\end{equation*}
Функцію $u$ називаємо (сильним) узагальненим розв'язком крайової задачі \eqref{1f1}, \eqref{1f2} з правою частиною $(f,g)$.

\textbf{Теорема 3.} \it Нехай параметри $\varphi\in\mathrm{RO}$ і $\lambda\in\mathbb{R}$ задовольняють нерівності $\sigma_0(\varphi)>m+1/2$ і $0<\lambda<\sigma_0(\varphi)-m+1/2$, а функції $\chi,\eta\in C^{\infty}(\overline{\Omega})$ задовольняють умову $\eta=1$ в околі $\mathrm{supp}\,\chi$. Тоді існує число $c=c(\varphi,\lambda,\chi,\eta)>0$ таке, що для довільної функції $u\in H^{\varphi}(\Omega)$ виконується оцінка
\begin{equation}\label{1f14}
\|\chi u\|_{H^{\varphi}(\Omega)}\leq c\,\bigl(\|\eta(A, B)u\|_{\mathcal{H}^{\varphi\varrho^{-2q}}(\Omega,
\Gamma)}+\|\eta u\|_{H^{\varphi\varrho^{-\lambda}}(\Omega)}\bigr).
\end{equation}
Тут $c$ не залежить від $u$. \rm

\textbf{\textit{Зауваження} 1.} У випадку, коли $\chi=\eta=1$, нерівність \eqref{1f14} є глобальною апріорною оцінкою узагальненого розв'язку $u$ еліптичної крайової задачі \eqref{1f1}, \eqref{1f2}. У цьому випадку умову $\lambda<\sigma_0(\varphi)-m+1/2$ можна прибрати.  Взагалі, нерівність \eqref{1f14} є локальною апріорною оцінкою розв'язку~$u$. Справді, для кожної непорожньої відкритої (у топології $\overline{\Omega}$) підмножини множини $\overline{\Omega}$, можна вибрати функції $\chi,\eta$ так, щоб вони задовольняли умову теореми~3 і їх носії лежали в цій підмножині. Якщо $0<\lambda\leq1$, то у нерівності \eqref{1f14} можна узяти $\chi(A, B)u$ замість $\eta(A, B)u$.

Дослідимо регулярність узагальнених розв'язків еліптичної крайової задачі \eqref{1f1}, \eqref{1f2}. Нехай $V$~--- відкрита множина в $\mathbb{R}^{n}$, яка має непорожній перетин з областю $\Omega$. Покладемо $\Omega_0:=\Omega\cap V$ і $\Gamma_{0}:=\Gamma\cap V$ (можливий випадок, коли $\Gamma_{0}=\varnothing$). Для довільного параметра $\alpha\in\mathrm{RO}$ введемо локальні аналоги просторів $H^{\alpha}(\Omega)$ і $H^{\alpha}(\Gamma)$.

За означенням, лінійний простір $H^{\alpha}_{\mathrm{loc}}(\Omega_{0},\Gamma_{0})$ cкладається з усіх розподілів $u\in\mathcal{D}'(\Omega)$ таких, що $\chi u\in H^{\alpha}(\Omega)$ для довільної функції $\chi\in C^{\infty}(\overline{\Omega})$ із $\mathrm{supp}\,\chi\subset\Omega_0\cup\Gamma_{0}$. Тут, як звичайно, $\mathcal{D}'(\Omega)$ позначає лінійний топологічний простір усіх розподілів в~$\Omega$. Топологія у лінійному просторі
$H^{\alpha}_{\mathrm{loc}}(\Omega_{0},\Gamma_{0})$ задається
напівнормами $u\mapsto\|\chi u\|_{H^{\alpha}(\Omega)}$, де $\chi$~--- довільна функція з означення цього простору.
Аналогічно, лінійний простір $H^{\alpha}_{\mathrm{loc}}(\Gamma_{0})$  складається з усіх розподілів $h\in\nobreak\mathcal{D}'(\Gamma)$ таких, що $\chi h\in H^{\alpha}(\Gamma)$ для довільної функції $\chi\in C^{\infty}(\Gamma)$ із $\mathrm{supp}\,\chi\subset\Gamma_{0}$.
Топологія у лінійному просторі
$H^{\alpha}_{\mathrm{loc}}(\Gamma_{0})$ задається
напівнормами $h\mapsto\|\chi h\|_{H^{\alpha}(\Gamma)}$, де $\chi$~--- довільна функція з означення цього простору.

\textbf{Теорема 4.} \it Нехай функція $u\in H^{m+1/2+}(\Omega)$ є узагальненим розв'язком еліптичної крайової задачі \eqref{1f1}, \eqref{1f2}, праві частини якої задовольняють умову
\begin{equation}\label{th4-cond}
(f,g)\in H^{\varphi\varrho^{-2q}}_{\mathrm{loc}}(\Omega_{0},\Gamma_{0})
\oplus\bigoplus_{j=1}^{q}
H^{\varphi\varrho^{-m_{j}-1/2}}_{\mathrm{loc}}(\Gamma_{0})=:
\mathcal{H}^{\varphi\varrho^{-2q}}_{\mathrm{loc}}(\Omega_{0},\Gamma_{0})
\end{equation}
для деякого функціонального параметра $\varphi\in\mathrm{RO}$ такого, що $\sigma_{0}(\varphi)>m+1/2$. Тоді розв'язок $u\in H^{\varphi}_{\mathrm{loc}}(\Omega_{0},\Gamma_{0})$.\rm

Відмітимо важливі окремі випадки цієї теореми. Якщо $\Omega_{0}=\Omega$ і $\Gamma_{0}=\Gamma$, то локальні простори $H^{\varphi}_{\mathrm{loc}}(\Omega_{0},\Gamma_{0})$ і $\mathcal{H}^{\varphi\varrho^{-2q}}_{\mathrm{loc}}(\Omega_{0},\Gamma_{0})$ збігаються з просторами $H^{\varphi}(\Omega)$ і $\mathcal{H}^{\varphi\varrho^{-2q}}(\Omega,\Gamma)$ відповідно. Тому теорема~4 стверджує, що регулярність узагальненого розв'язку $u$ підвищується глобально, тобто в усій області $\Omega$ аж до її межі $\Gamma$. Якщо  $\Gamma_{0}=\varnothing$ і $\Omega_{0}=\Omega$, то згідно з цією теоремою регулярність розв'язку $u$ підвищується в околах усіх внутрішніх точок замкненої області~$\overline\Omega$.

Теореми 1~-- 4 або їх версії відомі у випадку соболєвських просторів, коли $\varphi(\cdot)\equiv1$; див., наприклад, фундаментальну роботу С.~Агмона, А.~Дуглиса і Л.~Ніренберга \cite[розд.~5]{AgmonDouglisNirenberg59}, монографії Я.~А.~Ройтберга \cite[розд.~4, 7]{Roitberg96}, В.~О.~Козлова, В.~Г.~Маз'ї і Й.~Россмана \cite[розд.~4]{KozlovMazyaRossmann97}, Г.~Ескіна \cite[розд.~7]{Eskin11} та огляд М.~С.~Агарановича \cite[\S~2]{Agranovich97}. Відмітимо, що, мабуть, уперше Б.~Р.~Вайнберг і В.~В.~Грушин \cite[\S~4, формула~(76)]{VainbergGrushin67b} звернули увагу на те, що в описі \eqref{1f9} області значень оператора $(A,B)$ треба використовувати вираз вигляду
$$
\sum_{j=1}^{m-2q+1}(D_{\nu}^{j-1}f,w_{j})_{\Gamma}.
$$

Теореми 1~-- 4 і лему~1 доведемо у п.~7. Там же обґрунтуємо і зауваження~1.

\textbf{5. Застосування.} Як застосування просторів Хермандера дамо достатні умови неперервності узагальнених похідних (заданого порядку) розв'язків еліптичної крайової задачі \eqref{1f1}, \eqref{1f2}. Ці умови виводяться з теореми~4 і теореми вкладення Хермандера  \cite[теорема~2.2.7]{Hermander63}. Останню для розширеної соболєвської шкали можна сформулювати так: нехай $0\leq p\in\mathbb{Z}$ і  $\varphi\in\mathrm{RO}$, тоді
\begin{equation}\label{v}
\int\limits_1^{\infty} t^{2p+n-1}\varphi^{-2}(t)\,dt<\infty\;\;\Leftrightarrow\;\;
H^\varphi(\Omega)\subset C^p(\overline{\Omega}),
\end{equation}
причому вкладення неперервне; див. \cite[лема~2]{ZinchenkoMurach13UMJ11} або \cite[твердження~2.6(vi)]{MikhailetsMurach14}. Зауважимо, що у соболєвському випадку, коли $\varphi(t)\equiv t^{s}$ для деякого $s\in\mathbb{R}$, властивість \eqref{v} є теоремою вкладення Соболєва:
\begin{equation*}
s>p+n/2\;\;\Leftrightarrow\;\;
H^{(s)}(\Omega)\hookrightarrow C^p(\overline{\Omega}).
\end{equation*}

\textbf{Теорема 5.} \it Нехай ціле число $p\geq0$. Припустимо, що функція $u\in H^{m+1/2+}(\Omega)$ є узагальненим розв'язком еліптичної крайової задачі \eqref{1f1}, \eqref{1f2}, праві частини якої задовольняють умову \eqref{th4-cond} для деякого функціонального параметра $\varphi\in\mathrm{RO}$ такого, що $\sigma_{0}(\varphi)>m+1/2$~і
\begin{equation}\label{1f15}
\int\limits_1^{\infty}t^{2p+n-1}\varphi^{-2}(t)dt<\infty.
\end{equation}
Тоді $u\in C^{p}(\Omega_{0}\cup\Gamma_{0})$.\rm

\textbf{\textit{Зауваження} 2.} Умова \eqref{1f15} є точною у теоремі~5. А саме,  нехай $0\leq p\in\mathbb{Z}$, $\varphi\in\mathrm{RO}$ і $\sigma_{0}(\varphi)>m+1/2$; тоді з імплікації
\begin{equation}\label{implication}
\bigl(u\in H^{m+1/2+}(\Omega)\;\,\mbox{і}\;\,
(A,B)u\in\mathcal{H}^{\varphi\varrho^{-2q}}_
{\mathrm{loc}}(\Omega_{0},\Gamma_{0})\bigr)\;\Rightarrow\;
u\in C^{p}(\Omega_{0}\cup\Gamma_{0})
\end{equation}
випливає, що $\varphi$ задовольняє умову \eqref{1f15}.

Сформулюємо достатню умову, за якою узагальнений розв'язок $u$ крайової задачі \eqref{1f1}, \eqref{1f2} є класичним, тобто $u\in C^{2q}(\Omega)\cap C^{m}(U_{\sigma}\cup \Gamma)$ для деякого числа $\sigma>0$, де $U_{\sigma}:=\{x\in\Omega:\mathrm{dist}(x,\Gamma)<\sigma\}$. Якщо розв'язок $u$ цієї задачі класичний, то її ліві частини обчислюються за допомогою класичних похідних і є неперервними функціями на $\Omega$ і  $\Gamma$ відповідно.

\textbf{Теорема 6.} \it Нехай функція $u\in H^{m+1/2+}(\Omega)$ є узагальненим розв'язком еліптичної крайової задачі \eqref{1f1}, \eqref{1f2}, де
\begin{gather}\label{1f16}
f\in H^{\varphi_1\varrho^{-2q}}_{\mathrm{loc}}(\Omega,\varnothing)\cap H^{\varphi_2\varrho^{-2q}}_{\mathrm{loc}}(U_{\sigma},\Gamma),\\
g_j\in H^{\varphi_2\varrho^{-m_j-1/2}}(\Gamma)\quad\mbox{при кожному}\quad j\in\{1,\ldots,q\} \label{1f17}
\end{gather}
для деякого числа $\sigma>0$ і параметрів $\varphi_1,\varphi_2\in\mathrm{RO}$, які задовольняють умови $\sigma_0(\varphi_1)>m+1/2$, $\sigma_0(\varphi_2)>m+1/2$ і
\begin{gather}\label{f18}
\int\limits_1^{\infty}t^{2q+n-1}\varphi_{1}^{-2}(t)dt<\infty,\\
\int\limits_1^{\infty}t^{2m+n-1}\varphi_{2}^{-2}(t)dt<\infty. \label{f19}
\end{gather}
Тоді розв'язок $u$ класичний. \rm

Теореми 5, 6 і зауваження 2 будуть обґрунтовані у п.~7.

\textbf{6. Інтерполяція з функціональним параметром.} Простори Хермандера, які утворюють розширену соболєвську шкалу, можна отримати інтерполяцією з функціональним параметром пар гільбертових просторів Соболєва. Цей факт відіграватиме ключову роль у доведенні теореми~1. Метод інтерполяції з функціональним параметром гільбертових  просторів уперше з'явився у статті К.~Фояша і Ж.-Л.~Ліонса \cite[с.~278]{FoiasLions61}. Він є природнім узагальненням класичного інтерполяційного методу  Ж.-Л.~Ліонса \cite[розд.~1, п.~5]{LionsMagenes71} і С.-Г.~Крейна \cite[с.~253]{FunctionalAnalysis72} на випадок, коли  параметром інтерполяції служить не число, а досить загальна функція. Наведемо означення  інтерполяції з функціональним параметром пар гільбертових просторів та її властивості, потрібні у подальшому. Будемо слідувати монографії  \cite[п.~1.1]{MikhailetsMurach14}. Для наших цілей достатньо обмежитися випадком сепарабельних гільбертових просторів.

Нехай задана впорядкована пара $X:=[X_{0},X_{1}]$ сепарабельних комплексних гільбертових просторів $X_{0}$ і $X_{1}$ така, що $X_{1}$ є щільним лінійним многовидом у просторі $X_{0}$ та існує число $c>0$ таке, що $\|w\|_{X_{0}}\leq c\,\|w\|_{X_{1}}$ для довільного $w\in X_{1}$ (коротко кажучи, виконується неперервне і щільне вкладення $X_{1}\hookrightarrow X_{0}$). Пару $X$ називаємо припустимою. Для неї  існує самоспряжений додатно визначений оператор $J$ у гільбертовому просторі $X_{0}$ з областю визначення $X_{1}$ такий, що $\|Jw\|_{X_{0}}=\|w\|_{X_{1}}$ для довільного $w\in X_{1}$. Оператор $J$ називається породжуючим для $X$ і однозначно визначається за парою $X$.

Позначимо через $\mathcal{B}$ множину всіх вимірних за Борелем функцій
$\psi:\nobreak(0,\infty)\rightarrow(0,\infty)$, які відокремлені від нуля на кожній множині $[r,\infty)$ і обмежені на кожному відрізку $[a,b]$, де $r>0$ і $0<a<b<\infty$. Нехай $\psi\in\mathcal{B}$. У просторі $X_{0}$ за допомогою спектральної теореми означений, як функція від $J$, оператор $\psi(J)$, взагалі необмежений. Позначимо через $[X_{0},X_{1}]_\psi$ або, коротше, $X_{\psi}$ область визначення оператора $\psi(J)$, наділену скалярним добутком
$$
(w_1, w_2)_{X_\psi}:=(\psi(J)w_1,\psi(J)w_2)_{X_0}
$$
і відповідною нормою $\|w\|_{X_\psi}=(w,w)_{X_\psi}^{1/2}$. Простір $X_\psi$ гільбертів і сепарабельний, причому виконується неперервне і щільне вкладення $X_\psi \hookrightarrow X_0$.

Функцію $\psi\in\mathcal{B}$ називаємо інтерполяційним параметром, якщо для довільних припустимих пар $X=[X_0, X_1]$ і $Y=[Y_0, Y_1]$ гільбертових просторів і для будь-якого лінійного відображення $T$, заданого на $X_0$, виконується така умова: якщо при кожному $j\in\{0,1\}$ звуження відображення $T$ на простір $X_{j}$ є обмеженим оператором $T:X_{j}\rightarrow Y_{j}$, то і звуження відображення $T$ на простір $X_\psi$ є обмеженим оператором $T:X_{\psi}\rightarrow Y_{\psi}$. У~цьому випадку говоримо, що  простір $X_\psi$ отриманий інтерполяцією з функціональним параметром $\psi$ пари $X$.

Функція $\psi\in\mathcal{B}$ є інтерполяційним параметром тоді і тільки тоді, коли вона псевдоугнута в околі нескінченності, тобто еквівалентна там деякій угнутій додатній функції. Цей факт випливає з теореми Ж.~Петре \cite{Peetre68} про опис усіх інтерполяційних функцій додатного порядку.

Сформулюємо зазначену інтерполяційну властивість розширеної соболєвської шкали.

\textbf{Твердження 1.} \it Нехай задано функцію $\alpha\in\mathrm{RO}$ і дійсні числа $s_0$, $s_1$ такі, що $s_0<\sigma_0(\alpha)$ і $s_1>\sigma_1(\alpha)$. Покладемо
\begin{equation}\label{1f20}
\psi(t)=
\begin{cases}
\;t^{{-s_0}/{(s_1-s_0)}}\,
\alpha\bigl(t^{1/{(s_1-s_0)}}\bigr)&\text{при}\quad t\geq1, \\
\;\alpha(1)&\text{при}\quad0<t<1.
\end{cases}
\end{equation}
Тоді функція $\psi\in\mathcal{B}$ є інтерполяційним параметром і виконується така рівність просторів разом з еквівалентністю норм у них:
$$
\bigl[H^{(s_0)}(\Lambda),H^{(s_1)}(\Lambda)\bigr]_{\psi}=
H^{\alpha}(\Lambda),
$$
де $\Lambda\in\{\mathbb{R}^{n},\Omega,\Gamma\}$. Якщо $\Lambda=\mathbb{R}^{n}$, то буде рівність норм у цих просторах. \rm

Це твердження доведено в \cite[теореми 2.19 і 2.22]{MikhailetsMurach14} для $G\in\{\mathbb{R}^{n},\Gamma\}$ і в \cite[теорема~5.1]{MikhailetsMurach15ResMath1} для $G=\Omega$.

Відмітимо, що розширена соболєвська шкала замкнена відносно інтерполяції з функціональним параметром \cite[теорема~2.18]{MikhailetsMurach14} і збігається (з точністю до еквівалентності норм) з класом усіх гільбертових просторів, інтерполяційних для пар гільбертових просторів Соболєва \cite[теорема~2.24]{MikhailetsMurach14}. Остання властивість випливає з теореми В.~І.~Овчинникова \cite[п.~11.4]{Ovchinnikov84} про опис усіх гільбертових просторів, інтерполяційних для заданої пари гільбертових просторів. Нагадаємо, що властивість гільбертового простору $H$ бути інтерполяційним для припустимої пари $X=[X_0,X_1]$ значить таке: виконується неперервне вкладення $X_1\hookrightarrow H\hookrightarrow X_0$ і будь-який лініний оператор, обмеженим на кожному з просторів $X_0$ і $X_1$ є також обмеженим на~$H$.

Сформулюємо дві загальні властивості інтерполяції \cite[теореми 1.7, 1.5]{MikhailetsMurach14}, які будуть використані у наших доведеннях.

\textbf{Твердження 2.} \it Нехай задано дві припустимі пари $X=[X_0,X_1]$ і $Y=[Y_0,Y_1]$ гільбертових просторів. Нехай, окрім того, на $X_0$ задано лінійне відображення $T$ таке, що його звуження на простори $X_j$, де $j=0,1$, є обмеженими і нетеровими операторами $T:X_j\rightarrow Y_j$, які мають спільне ядро і однаковий індекс. Тоді для довільного інтерполяційного параметра $\psi\in\mathcal{B}$ обмежений оператор $T:X_\psi\rightarrow Y_\psi$ нетерів з тим же ядром і індексом, а його область значень дорівнює $Y_\psi\cap T(X_0)$. \rm

\textbf{Твердження 3.} \it Нехай задано скінченне число припустимих пар
$[X_{0}^{(j)},X_{1}^{(j)}]$  гільбертових просторів, де $j=1,\ldots,q$. Тоді для довільної функції $\psi\in\mathcal{B}$ правильна така рівність просторів разом з рівністю норм у них:
\begin{equation*}
\biggl[\,\bigoplus_{j=1}^{q}X_{0}^{(j)},\,
\bigoplus_{j=1}^{q}X_{1}^{(j)}\biggr]_{\psi}=\,
\bigoplus_{j=1}^{q}\bigl[X_{0}^{(j)},\,X_{1}^{(j)}\bigr]_{\psi}.
\end{equation*} \rm

Лінійні диференціальні оператори з гладкими коефіцієнтами є обмеженими на парах підходящих просторів Хермандера. А саме, є правильним такий результат.

\textbf{Твердження 4.} \it $\mathrm{(i)}$ Нехай $L$ є лінійний диференціальний вираз порядку $l\geq0$ на $\overline{\Omega}$ з коефіцієнтами класу $C^{\infty}(\overline{\Omega})$. Тоді відображення $u\mapsto Lu$, де $u\in C^{\infty}(\overline{\Omega})$, продовжується єдиним чином (за неперервністю) до обмеженого лінійного оператора
\begin{equation*}
L:H^\alpha(\Omega)\rightarrow H^{\alpha\varrho^{-l}}(\Omega)
\end{equation*}
для кожного параметра $\alpha\in\mathrm{RO}$.

$\mathrm{(ii)}$ Нехай $K$ є крайовий лінійний диференціальний вираз порядку $k\geq\nobreak0$ на межі $\Gamma$ з коефіцієнтами класу $C^{\infty}(\Gamma)$. Тоді відображення $u\mapsto Ku$, де $u\in C^{\infty}(\overline{\Omega})$, продовжується єдиним чином (за неперервністю) до обмеженого лінійного оператора
\begin{equation*}
K:H^\alpha(\Omega)\rightarrow H^{\alpha\varrho^{-k-1/2}}(\Gamma)
\end{equation*}
для кожного параметра $\alpha\in\mathrm{RO}$ такого, що $\sigma_{0}(\alpha)>k+1/2$. \rm

У випадку просторів Соболєва, коли $\alpha(t)\equiv t^{s}$, твердження 4 добре відоме. Звідси випадок довільного $\alpha\in\mathrm{RO}$  виводиться за допомогою інтерполяції з функціональним параметром на підставі твердження 1.

\textbf{7. Доведення.} Доведемо теореми 1~-- 6, лему 1 та обґрунтуємо зауваження 1 і~2.

\textbf{\textit{Доведення теореми} 1.} У соболєвському випадку, коли $\varphi=\varrho^{s}$ і дійсне $s>m+1/2$, ця теорема відома за виключенням вказаного зв'язку скінченновимірного простору $N_{\star}$ з формально спряженою задачею \eqref{1f4}, \eqref{1f5}. У такому вигляді теорема~1 міститься у результаті, доведеному в монографії Я.~А.~Ройтберга \cite[теорема~4.1.3]{Roitberg96}. У повному обсязі, але за додаткового припущення $s\in\mathbb{Z}$, теорема~1 міститься у результаті, встановленому в монографії В.~О.~Козлова, В.~Г.~Маз'ї і Й.~Россмана \cite[наслідок 4.1.1]{KozlovMazyaRossmann97}. Покажемо, що і для дробових $s$ висновок цієї теореми правильний у повному обсязі.

Згідно з \cite[теорема~4.1.3]{Roitberg96} відображення \eqref{1f3} продовжується за неперервністю до обмеженого і нетерового оператора
\begin{equation}\label{roit-oper}
(A,B):H^{s,(m+1)}(\Omega)\rightarrow H^{s-2q,(m+1-2q)}(\Omega) \oplus\bigoplus_{j=1}^{q}H^{(s-m_j-1/2)}(\Gamma)=:
\mathcal{Q}^{s-2q}(\Omega,\Gamma)
\end{equation}
для довільного $s\in\mathbb{R}$. Тут $H^{s,(r)}(\Omega)$, де $s\in\mathbb{R}$ і $1\leq r\in\mathbb{Z}$, є модифікований за Ройтбергом гільбертів простір Соболєва \cite[п.~2.1]{Roitberg96}. Зокрема, якщо $s\geq0$ і $s\notin\{1/2,\ldots,r-1/2\}$, то $H^{s,(r)}(\Omega)$ є, за означенням, поповненням простору $C^{\infty}(\overline{\Omega})$ за нормою
$$
\|u\|_{H^{s,(r)}(\Omega)}:=
\biggl(\|u\|_{H^{(s)}(\Omega)}^{2}+
\sum_{k=1}^{r}\;\|(D_{\nu}^{k-1}u)\!\upharpoonright\!\Gamma\|
_{H^{(s-k+1/2)}(\Gamma)}^{2}\biggr)^{1/2}.
$$
Відмітимо, що виконується неперервне вкладення $H^{s+\delta,(r)}(\Omega)\hookrightarrow H^{s,(r)}(\Omega)$ при $\delta>0$. Окрім того, якщо $s>r-1/2$, то простори $H^{s,(r)}(\Omega)$ і $H^{(s)}(\Omega)$ рівні як поповнення $C^{\infty}(\overline{\Omega})$ за еквівалентними нормами. Тому оператор \eqref{1f8}, де $\varphi=\varrho^{s}$, і оператор \eqref{roit-oper} рівні при $s>r-1/2$.

Згідно зі згаданим результатом \cite[теорема~4.1.3]{Roitberg96} ядро оператора  \eqref{roit-oper} збігається з $N$, а область значень складається з усіх векторів $(f,g_1,\ldots,g_q)\in\mathcal{Q}^{s-2q}(\Omega,\Gamma)$ таких, що задовольняють умову \eqref{1f9}, у якій замість $N_{\star}$ фігурує
деякий скінченновимірний простір, що лежить в $C^{\infty}(\overline{\Omega})\times(C^{\infty}(\Gamma))^{m-q+1}$ і не залежить від $s$. Звідси негайно випливає рівність
\begin{equation*}
(A,B)(H^{s_{2},(m+1)}(\Omega))=\mathcal{Q}^{s_{2}-2q}(\Omega,\Gamma)\cap
(A,B)(H^{s_{1},(m+1)}(\Omega))\quad\mbox{при}\quad s_{1}<s_{2}.
\end{equation*}
Зокрема,
\begin{equation}\label{proof-th1-a}
(A,B)(H^{(s)}(\Omega))=\mathcal{Q}^{s-2q}(\Omega,\Gamma)\cap
(A,B)(H^{m,(m+1)}(\Omega))\quad\mbox{при}\quad m+1/2<s\in\mathbb{R}.
\end{equation}
Згідно з \cite[теорема 4.1.4]{KozlovMazyaRossmann97} простір $(A,B)(H^{m,(m+1)}(\Omega))$ складається з усіх векторів $(f,g_1,\ldots,g_q)\in\mathcal{Q}^{m-2q}(\Omega,\Gamma)$, які задовольняють умову \eqref{1f9}, де $(\cdot,\cdot)_{\Gamma}$ є також продовженням за неперервністю скалярного добутку в $L_{2}(\Gamma)$. Тому для кожного дійсного $s>m+1/2$ область значень $(A,B)(H^{(s)}(\Omega))$ оператора \eqref{1f8}, де $\varphi=\varrho^{s}$, є такою як це стверджується у теоремі~1. Отже, у соболєвському випадку ця теорема обґрунтована.

У загальній ситуації доведемо її за допомогою інтерполяції з функціональним параметром пар деяких просторів Соболєва. За умовою, $\varphi\in\mathrm{RO}$ та $\sigma_0(\varphi)>m+1/2$. Виберемо дійсні числа $l_{0}$ і $l_{1}$ такі, що
$m+1/2<l_{0}<\sigma_{0}(\varphi)$ і $\sigma_{1}(\varphi)<l_{1}$.
Відображення \eqref{1f3} продовжується за неперервністю до
обмежених і нетерових операторів
\begin{equation}\label{01f21}
(A,B):\,H^{(l_{i})}(\Omega)\rightarrow
H^{(l_{i}-2q)}(\Omega)\oplus
\bigoplus_{j=1}^{q} H^{(l_{i}-m_{j}-1/2)}(\Gamma)=:
\mathcal{H}^{(l_{i}-2q)}(\Omega,\Gamma),\quad
i\in\{0,1\},
\end{equation}
які діють у соболєвських просторах. Ці оператори мають спільне ядро $N$ і однаковий індекс, який дорівнює $\dim N-\dim N_{\star}$. Окрім того,
\begin{gather}\label{01f22}
(A,B)(H^{(l_{i})}(\Omega))=
\bigl\{(f,g)\in\mathcal{H}^{(l_{i}-2q)}(\Omega,\Gamma):
\mbox{виконується \eqref{1f9}}\bigr\}.
\end{gather}

Означимо функцію $\psi\in\mathcal{B}$ за формулою \eqref{1f20}, у якій покладемо $\alpha:=\varphi$. Ця функція є інтерполяційним
параметром згідно з твердженням~1. Тому на підставі твердження~2 з обмеженості та нетеровості обох операторів \eqref{01f21} випливає обмеженість і нетеровість оператора
\begin{equation}\label{01f23}
(A,B):\bigl[H^{(l_{0})}(\Omega),H^{(l_{1})}(\Omega)\bigr]_{\psi}\to
\bigl[\mathcal{H}^{(l_{0}-2q)}(\Omega,\Gamma),
\mathcal{H}^{(l_{1}-2q)}(\Omega,\Gamma)\bigr]_{\psi}.
\end{equation}
Він є звуженням оператора \eqref{01f21} з $i=0$. Покажемо, що
\eqref{01f23}~--- це оператор \eqref{1f8} із теореми~1.

На підставі тверджень 1 і 3 маємо такі рівності просторів разом з еквівалентністю норм у них:
\begin{gather*}
\bigl[H^{(l_{0})}(\Omega),H^{(l_{1})}(\Omega)\bigr]_{\psi}=
H^{\varphi}(\Omega),\\
\bigl[\mathcal{H}^{(l_{0}-2q)}(\Omega,\Gamma),
\mathcal{H}^{(l_{1}-2q)}(\Omega,\Gamma)\bigr]_{\psi}=
\bigl[H^{(l_{0}-2q)}(\Omega),H^{(l_{1}-2q)}(\Omega)\bigr]_{\psi}\oplus\\
\oplus\bigoplus_{j=1}^{q}\bigl[H^{(l_{0}-m_{j}-1/2)}(\Gamma),
H^{(l_{1}-m_{j}-1/2)}(\Gamma)\bigr]_{\psi}=
\mathcal{H}^{\varphi\varrho^{-2q}}(\Omega,\Gamma).
\end{gather*}
Тому обмежений і нетерів оператор \eqref{01f23} діє в парі просторів \eqref{1f8}. Оскільки цей оператор є продовженням за неперервністю відображення \eqref{1f3}, то він є оператором \eqref{1f8}. На підставі твердження~2 ядро цього оператора та його індекс збігаються з спільним ядром $N$ і однаковим індексом $\dim N-\dim N_{\star}$ операторів \eqref{01f21}. Окрім того, область значень оператора \eqref{1f8} дорівнює
\begin{gather*}
(A,B)(H^{\varphi}(\Omega))=
\mathcal{H}^{\varphi\varrho^{-2q}}(\Omega,\Gamma)\cap
(A,B)(H^{(l_{0})}(\Omega))=\\
=\bigl\{(f,g)\in\mathcal{H}^{\varphi\varrho^{-2q}}(\Omega,\Gamma):\,
\mbox{виконується}\;\eqref{1f9}\bigr\}.
\end{gather*}
Тут використали рівність \eqref{01f22} та вкладення
$$
\mathcal{H}^{\varphi\varrho^{-2q}}(\Omega,\Gamma)\hookrightarrow
\mathcal{H}^{(l_{0}-2q)}(\Omega,\Gamma),
$$
яке випливає з властивості~\eqref{1f7}, оскільки $l_{0}<\sigma_{0}(\varphi)$. Таким чином, доведено всі властивості оператора \eqref{1f8}, сформульовані в теоремі 1.

Теорема 1 доведена.

\textbf{\textit{Доведення леми} 1.} Скористаємося обмеженим нетеровим оператором \eqref{roit-oper} для $s:=m$. Згідно з \cite[теорема 4.1.4]{KozlovMazyaRossmann97} вимірність коядра цього оператора дорівнює $\dim N_{\star}$. Лінійний многовид  $C^{\infty}(\overline{\Omega})\times(C^{\infty}(\Gamma))^{q}$ щільний у просторі $\mathcal{Q}^{m-2q}(\Omega,\Gamma)$.
Тому на підставі \cite[лема 2.1]{HohbergKrein57} існує скінченновимірний простір $G\subset C^{\infty}(\overline{\Omega})\times(C^{\infty}(\Gamma))^{q}$ такий, що
\begin{equation}\label{proof-lemma-a}
\mathcal{Q}^{m-2q}(\Omega,\Gamma)=G\dotplus (A,B)(H^{m,(m+1)}(\Omega)).
\end{equation}
Звідси випливає, що $\dim G=\dim N_{\star}$.

Нехай число $s$ задовольняє умову $m+1/2<s<\sigma_{0}(\varphi)$. Тоді виконується неперервні вкладення
$$
\mathcal{H}^{\varphi\varrho^{-2q}}(\Omega,\Gamma)\hookrightarrow
\mathcal{H}^{(s-2q)}(\Omega,\Gamma)=\mathcal{Q}^{s-2q}(\Omega,\Gamma)
\hookrightarrow\mathcal{Q}^{m-2q}(\Omega,\Gamma)
$$
на підставі \eqref{1f7} і того, що простори $H^{s-2q,(m+1-2q)}(\Omega)$ і $H^{(s-2q)}(\Omega)$ рівні з точністю до еквівалентності норм при $s-2q>m+1-2q-1/2$, як це зазначалося у доведенні теореми~1. Окрім того, $G\subset\mathcal{H}^{\varphi\varrho^{-2q}}(\Omega,\Gamma)$. Тому з рівності \eqref{proof-lemma-a} випливає формула
\begin{equation}\label{proof-lemma-b}
\mathcal{H}^{\varphi\varrho^{-2q}}(\Omega,\Gamma)=G\dotplus
\bigl((A,B)(H^{m,(m+1)}(\Omega))\cap
\mathcal{H}^{\varphi\varrho^{-2q}}(\Omega,\Gamma)\bigr).
\end{equation}
Згідно з \cite[теорема 4.1.4]{KozlovMazyaRossmann97} область значень  $(A,B)(H^{m,(m+1)}(\Omega)$ оператора \eqref{roit-oper}, де $s=m$, складається з усіх векторів $(f,g_1,\ldots,g_q)\in\mathcal{Q}^{m-2q}(\Omega,\Gamma)$, які задовольняють умову \eqref{1f9}. Тому другий доданок у сумі \eqref{proof-lemma-b} складається з усіх векторів $(f,g_1,\ldots,g_q)\in\mathcal{H}^{\varphi\varrho^{-2q}}(\Omega,\Gamma)$, які задовольняють \eqref{1f9}. Отже, \eqref{proof-lemma-b} перетворюється на рівність \eqref{1f11}. У~ній, згідно з нашими міркуваннями, простір $G$ не залежить від $s$.

Лема 1 доведена.

\textbf{\textit{Доведення теореми}~2.} Згідно з теоремою~1, $N$~--- ядро, а $Q(\mathcal{H}^{\varphi\varrho^{-2q}}(\Omega,\Gamma))$~--- область значень оператора \eqref{1f8}. Тому звуження відображення \eqref{1f8} на простір $P(H^{\varphi}(\Omega))$ є обмеженим лінійним бієктивним оператором. Отже, він є ізоморфізмом \eqref{1f12} за теоремою Банаха про обернений оператор. Теорема 2 доведена.

\textbf{\textit{Доведення теореми}~3.} У випадку, коли $\chi=\eta=1$, ця теорема є наслідком скінченновимірності ядра і замкненості області значень оператора \eqref{1f8}, доведених у теоремі~1, та компактності вкладення $H^{\varphi\varrho^{-\lambda}}(\Omega)\hookrightarrow H^{\varphi}(\Omega)$. Це стверджує лема Пітре \cite[лема~3]{Peetre61}. У цьому випадку $\lambda$~--- довільне додатне число. Таким чином, існує число $\tilde{c}=\tilde{c}(\varphi,\lambda)>0$ таке, що для довільної функції $v\in H^{\varphi}(\Omega)$ виконується глобальна апріорна оцінка
\begin{equation}\label{global-estimate}
\|v\|_{H^{\varphi}(\Omega)}\leq\tilde{c}\,
\bigl(\|(A, B)v\|_{\mathcal{H}^{\varphi\varrho^{-2q}}(\Omega,\Gamma)}+
\|v\|_{H^{\varphi\varrho^{-\lambda}}(\Omega)}\bigr).
\end{equation}

Виведемо з цієї оцінки теорему~3 для $\lambda=1$. Зауважимо спочатку, що нерівність $\lambda<\sigma_0(\varphi)-m+1/2$, вказана в умові цієї теореми, виконується для $\lambda=1$. Довільно виберемо функцію $u\in H^{\varphi}(\Omega)$. Нехай функції $\chi,\eta\in C^{\infty}(\overline{\Omega})$ такі як в умові теореми~3. Узявши $v:=\chi u\in H^{\varphi}(\Omega)$ і $\lambda:=1$ в оцінці \eqref{global-estimate}, запишемо
\begin{equation}\label{1f29}
\|\chi u\|_{H^{\varphi}(\Omega)}\leq\tilde{c}\,
\bigl(\|(A, B)(\chi u)\|_{\mathcal{H}^{\varphi\varrho^{-2q}}(\Omega,
\Gamma)}+\|\chi u\|_{H^{\varphi\varrho^{-1}}(\Omega)}\bigr).
\end{equation}
Переставивши оператор множення на функцію $\chi$ з
диференціальними операторами $A$ і $B_{1},\ldots,B_{q}$, отримаємо рівність
\begin{equation*}
(A,B)(\chi u)=(A,B)(\chi\eta u)=\chi(A,B)(\eta u)+(A',B')(\eta u)=
\chi(A,B)u+(A',B')(\eta u).
\end{equation*}
Тут $A'$~--- деякий лінійний диференціальний оператор на $\overline{\Omega}$ порядку $\mathrm{ord}\,A'\leq 2q-1$, а $B':=(B_{1}',\ldots,B_{q}')$~--- набір деяких крайових лінійних диференціальних операторів на $\Gamma$, порядки яких задовольняють умову $\mathrm{ord}\,B_{j}'\leq m_{j}-1$ для кожного $j\in\{1,\ldots,q\}$. При цьому всі коефіцієнти операторів $A'$ і $B_{j}'$ належать до $C^{\infty}(\overline{\Omega})$ і $C^{\infty}(\Gamma)$ відповідно. Таким чином,
\begin{equation}\label{1f30}
(A,B)(\chi u)=\chi(A,B)u+(A',B')(\eta u).
\end{equation}
Згідно з твердженням 4 виконується нерівність
\begin{equation}\label{1f30b}
\|(A',B')(\eta u)\|_{\mathcal{H}^{\varphi\varrho^{-2q}}(\Omega,\Gamma)}
\leq c_{1}\|\eta u\|_{H^{\varphi\varrho^{-1}}(\Omega)}.
\end{equation}
Тут і далі у доведенні через $c_{1},\ldots,c_{7}$ позначено деякі додатні числа, не залежні від~$u$.

На підставі формул \eqref{1f29}~-- \eqref{1f30b} отримаємо нерівності
\begin{gather*}
\|\chi u\|_{H^{\varphi}(\Omega)}\leq \tilde{c}\,
\bigl(\|\chi(A,B)u\|_{\mathcal{H}^{\varphi\varrho^{-2q}}(\Omega,\Gamma)}+
\|(A',B')(\eta u)\|_{\mathcal{H}^{\varphi\varrho^{-2q}}(\Omega,\Gamma)}
+\|\chi u\|_{H^{\varphi\varrho^{-1}}(\Omega)}\bigr)\leq\\
\leq\tilde{c}\,
\|\chi(A, B)u\|_{\mathcal{H}^{\varphi\varrho^{-2q}}(\Omega,\Gamma)}+
\tilde{c}\,c_{1}\|\eta u\|_{H^{\varphi\varrho^{-1}}(\Omega)}+
\tilde{c}\,\|\chi u\|_{H^{\varphi\varrho^{-1}}(\Omega)}.
\end{gather*}
Тут на підставі твердження 4
\begin{gather*}
\|\chi u\|_{H^{\varphi\varrho^{-1}}(\Omega)}=
\|\chi\eta u\|_{H^{\varphi\varrho^{-1}}(\Omega)}\leq
c_{2}\|\eta u\|_{H^{\varphi\varrho^{-1}}(\Omega)}.
\end{gather*}
Отже,
\begin{equation}\label{proof-th3}
\|\chi u\|_{H^{\varphi}(\Omega)}\leq c_{3}
\bigl(\|\chi(A,B)u\|_{\mathcal{H}^{\varphi\varrho^{-2q}}(\Omega,\Gamma)}+
\|\eta u\|_{H^{\varphi\varrho^{-1}}(\Omega)}).
\end{equation}
З цієї нерівності випливає потрібна оцінка \eqref{1f14}, оскільки
\begin{equation*}
\|\chi(A,B)u\|_{\mathcal{H}^{\varphi\varrho^{-2q}}(\Omega,\Gamma)}=
\|\chi\eta(A,B)u\|_{\mathcal{H}^{\varphi\varrho^{-2q}}(\Omega,\Gamma)}\leq
c_{4}\|\eta(A,B)u\|_{\mathcal{H}^{\varphi\varrho^{-2q}}(\Omega,\Gamma)}
\end{equation*}
на підставі твердження 4. Теорема 3 доведена у випадку, коли $\lambda=1$. Звісно, її висновок правильний і якщо $0<\lambda<1$.

Доведемо тепер цю теорему у випадку, коли
\begin{equation}\label{proof-th3-a}
1<\lambda<\sigma_0(\varphi)-m+1/2.
\end{equation}
Для кожного дійсного числа $l\geq1$ позначимо через $\mathcal{P}_{l}$ висновок теореми~3 у випадку, коли $\lambda=l$. А саме, $\mathcal{P}_{l}$ позначає таке твердження: для довільних функцій $\varphi\in\mathrm{RO}$ і $\chi,\eta\in C^{\infty}(\overline{\Omega})$, які задовольняють умови $\sigma_0(\varphi)>m+1/2$, $l<\sigma_0(\varphi)-m+1/2$ і $\eta=1$ в околі $\mathrm{supp}\,\chi$, існує число $c=c(\varphi,l,\chi,\eta)>0$ таке, що для довільної функції $u\in H^{\varphi}(\Omega)$ виконується нерівність \eqref{1f14} з $\lambda=l$. Істинність твердження $\mathcal{P}_{1}$ доведена вище. Довільно виберемо дійсні числа $l\geq1$ і $\delta\in(0,1]$. Доведемо, що $\mathcal{P}_{l}\Rightarrow\mathcal{P}_{l+\delta}$.

Припустимо, що твердження $\mathcal{P}_{l}$ істинне. Нехай функції $\varphi\in\mathrm{RO}$ і $\chi,\eta\in C^{\infty}(\overline{\Omega})$ задовольняють умови $\sigma_0(\varphi)>m+1/2$, $l+\delta<\sigma_0(\varphi)-m+1/2$ і $\eta=1$ в околі $\mathrm{supp}\,\chi$. Тоді знайдеться функція $\eta_{1}\in C^{\infty}(\overline{\Omega})$ така, що $\eta_{1}=1$ в околі $\mathrm{supp}\,\chi$ і $\eta=1$ в околі $\mathrm{supp}\,\eta_{1}$. За припущенням, існує число $c_{5}>0$ таке, що для довільної функції $u\in H^{\varphi}(\Omega)$ виконується оцінка
\begin{equation}\label{proof-th3-b}
\|\chi u\|_{H^{\varphi}(\Omega)}\leq c_{5}\bigl
(\|\eta_{1}(A,B)u\|_{\mathcal{H}^{\varphi\varrho^{-2q}}(\Omega,\Gamma)}
+\|\eta_{1}u\|_{H^{\varphi\varrho^{-l}}(\Omega)}\bigr).
\end{equation}
Оскільки $\sigma_0(\varphi\varrho^{-l-\delta+1})>m+1/2$, то на підставі твердження $\mathcal{P}_{1}$ маємо оцінку
\begin{equation}\label{proof-th3-c}
\begin{gathered}
\|\eta_{1}u\|_{H^{\varphi\varrho^{-l}}(\Omega)}\leq
\|\eta_{1}u\|_{H^{\varphi\varrho^{-l-\delta+1}}(\Omega)}\leq\\
\leq c_{6}\bigl(\|\eta(A,B)u\|_
{\mathcal{H}^{\varphi\varrho^{-l-\delta+1-2q}}(\Omega,\Gamma)}
+\|\eta u\|_{H^{\varphi\varrho^{-l-\delta}}(\Omega)}\bigr).
\end{gathered}
\end{equation}
Окрім того,
\begin{equation}\label{proof-th3-d}
\|\eta_{1}(A,B)u\|_{\mathcal{H}^{\varphi\varrho^{-2q}}(\Omega,\Gamma)}=
\|\eta_{1}\eta(A,B)u\|_{\mathcal{H}^{\varphi\varrho^{-2q}}(\Omega,\Gamma)}
\leq c_{7}\|\eta(A,B)u\|_{\mathcal{H}^{\varphi\varrho^{-2q}}(\Omega,\Gamma)}.
\end{equation}
На підставі оцінок \eqref{proof-th3-b}~-- \eqref{proof-th3-d} запишемо
\begin{equation*}
\|\chi u\|_{H^{\varphi}(\Omega)}\leq c_{5}c_{7}
\|\eta(A,B)u\|_{\mathcal{H}^{\varphi\varrho^{-2q}}(\Omega,\Gamma)}+
c_{5}c_{6}\bigl(\|\eta(A,B)u\|_
{\mathcal{H}^{\varphi\varrho^{-2q}}(\Omega,\Gamma)}
+\|\eta u\|_{H^{\varphi\varrho^{-l-\delta}}(\Omega)}\bigr),
\end{equation*}
тобто отримали нерівність \eqref{1f14} з $\lambda=l+\delta$. Імплікація  $\mathcal{P}_{l}\Rightarrow\mathcal{P}_{l+\delta}$ обґрунтована.

Тепер можемо довести теорему~3 у випадку \eqref{proof-th3-a}. За доведеним, правильний ланцюжок імплікацій
\begin{equation*}
\mathcal{P}_{1}\Rightarrow\mathcal{P}_{2}\Rightarrow\ldots
\Rightarrow\mathcal{P}_{[\lambda]}\Rightarrow \mathcal{P}_{\lambda},
\end{equation*}
де твердження $\mathcal{P}_{1}$ істинне, а $\mathcal{P}_{\lambda}$ є висновком теореми~3 у досліджуваному випадку (як звичайно, $[\lambda]$~--- ціла частина числа $\lambda$). Тому цей висновок є також істинним.

Теорема~3 доведена.

У зауваженні~1 потребують обґрунтування друге і останнє речення. Друге речення обґрунтоване у першому абзаці доведення цієї теореми, а останнє речення є прямим наслідком оцінки \eqref{proof-th3}.

\textbf{\textit{Доведення теореми} 4.} Спочатку обґрунтуємо цю теорему у випадку, коли $\Omega_{0}=\Omega$ і $\Gamma_{0}=\Gamma$. За умовою,
$u\in H^{(s)}(\Omega)$ для деякого дійсного числа $s$ такого, шо $m+1/2<s<\sigma_{0}(\varphi)$, і $(f,g)=(A,B)u\in\mathcal{H}^{\varphi\rho^{-2q}}(\Omega,\Gamma)$. Тому
$$
(f,g)\in\mathcal{H}^{\varphi\varrho^{-2q}}(\Omega,\Gamma)
\cap(A,B)(H^{(s)}(\Omega))=(A,B)(H^{\varphi}(\Omega));
$$
тут рівність правильна на підставі теореми~1. Отже, поряд з умовою $(A,B)u=(f,g)$ виконується рівність $(A,B)v=(f,g)$ для деякого
$v\in H^{\varphi}(\Omega)$. Тому $(A,B)(u-v)=0$, що за теоремою 1 тягне за собою включення $w:=u-v\in N\subset C^{\infty}(\overline{\Omega})$. Звідси $u=v+w\in H^{\varphi}(\Omega)$. У досліджуваному випадку теорема~4 доведена.

Доведемо її в загальному випадку. Міркування проведемо за схемою, наведеною в \cite[с.~308]{AnopKasirenko16MFAT4}. Довільно виберемо відкриту множину $V_{1}\subset\mathbb{R}^{n}$ таку, що $\overline{V_{1}}\subset V$ і $\Omega\cap V_{1}\neq\varnothing$ та покладемо $\Omega_1:=\Omega\cap V_1$ і $\Gamma_{1}:=\Gamma\cap V_1$. Доведемо, що $u\in H^{\varphi}_{\mathrm{loc}}(\Omega_{1},\Gamma_{1})$.

Нехай функції $\chi,\eta\in C^{\infty}(\overline{\Omega})$ такі, що їх носії лежать в $\Omega_{0}\cup\Gamma_{0}$ і $\chi=1$ в околі  $\Omega_{1}\cup\Gamma_{1}$ та $\eta=1$ на $\mathrm{supp}\,\chi$. За умовою, $u\in H^{(s)}(\Omega)$ для деякого $s\in\mathbb{R}$ такого, що $m+1/2<s<\sigma_{0}(\varphi)$, і $(A,B)u=(f,g)\in \mathcal{H}^{\varphi\rho^{-2q}}_{\mathrm{loc}}(\Omega_{0},\Gamma_{0})$. Тому
\begin{equation*}
(A,B)(\chi u)=\eta(A,B)(\chi u)=\eta(f,g)-\eta(A,B)((1-\chi)u).
\end{equation*}
Використовуючи проектор $P_{\star}$ з теореми~2, запишемо $(A,B)(\chi u)=P_{\star}(\eta(f,g))+F$, де
\begin{equation*}
F:=(1-P_{\star})(\eta(f,g))-\eta(A,B)((1-\chi)u).
\end{equation*}
Оскільки $P_{\star}(\eta(f,g))\in P_{\star}(\mathcal{H}^{\varphi\rho^{-2q}}(\Omega,\Gamma))$,
то
\begin{equation*}
F=(A,B)(\chi u)-P_{\star}(\eta(f,g))\in
P_{\star}\bigl(\mathcal{H}^{\varrho^{s-2q}}(\Omega,\Gamma)\bigr).
\end{equation*}
За теоремою~2, існують функції $u_1\in H^{\varphi}(\Omega)$ і $u_2\in H^{(s)}(\Omega)$ такі, що $(A,B)u_1=P_{\star}(\eta(f,g))$ і $(A,B)u_2=F$. Тоді $(A,B)(\chi u-u_1-u_2)=0$, звідки
\begin{equation*}
w:=\chi u-u_1-u_2\in N\subset C^{\infty}(\overline{\Omega})
\end{equation*}
на підставі теореми~1. Помітимо, що  $F\in\mathcal{H}^{\varrho^{l-2q}}_{\mathrm{loc}}(\Omega_{1},\Gamma_{1})$ для кожного дійсного числа $l>\sigma_{1}(\varphi)$, оскільки $(1-P_{\star})(\eta(f,g))\in N_{\star}$ і $\eta(A,B)((1-\chi)u)=0$ на $\Omega_{1}\cup\Gamma_{1}$. Тому
\begin{equation*}
u_2\in H_{\mathrm{loc}}^{\varrho^l}(\Omega_{1},\Gamma_{1})\subset
H^{\varphi}_{\mathrm{loc}}(\Omega_{1},\Gamma_{1})
\end{equation*}
згідно з теоремою про локальне підвищення регулярності розв'язків еліптичних крайових задач у просторах Соболєва (див., наприклад, \cite[теорема 7.2.1]{Roitberg96}). Таким чином,
\begin{equation*}
\chi u=u_1+u_2+w\in H^{\varphi}_{\mathrm{loc}}(\Omega_{1},\Gamma_{1}).
\end{equation*}
Отже, $\zeta u=\zeta\chi u\in H^{\varphi}(\Omega)$ для довільної функції $\zeta\in C^{\infty}(\overline{\Omega})$, яка задовольняє умову  $\mathrm{supp}\,\zeta\subset\Omega_{1}\cup\Gamma_{1}$, тобто $u\in H^{\varphi}_{\mathrm{loc}}(\Omega_{1},\Gamma_{1})$. Тепер $u\in H^{\varphi}_{\mathrm{loc}}(\Omega_{0},\Gamma_{0})$ згідно із зробленим вибором множини $V_{1}$.

Теорема 4 доведена.

\textbf{\textit{Доведення теореми}~5.} Довільно виберемо точку
$x\in\Omega_{0}\cup\Gamma_{0}$ і функцію $\chi\in C^{\infty}(\overline{\Omega})$ таку, що $\mathrm{supp}\,\chi\subset\Omega_0\cup\Gamma_0$ і $\chi=1$ у деякому околі $V(x)$ точки~$x$. З теореми 4, умови \eqref{1f15} і еквівалентності~\eqref{v} випливає включення $\chi u\in H^{\varphi}(\Omega)\subset C^{p}(\overline{\Omega})$.
Тому $u\in C^{p}(V(x))$. Звідси, з урахуванням довільності вибору точки $x$, робимо висновок, що $u\in C^{p}(\Omega_{0}\cup\Gamma_{0})$.
Теорема~5 доведена.

Обґрунтуємо зауваження~2. Нехай $0\leq p\in\mathbb{Z}$, $\varphi\in\mathrm{RO}$ і $\sigma_{0}(\varphi)>m+1/2$. Припустимо, що імплікація \eqref{implication} істинна. Нехай $V$~--- деяка відкрита куля така, що  $\overline{V}\subset\Omega_{0}$. Довільно виберемо функцію $v\in H^{\varphi}(V)$. Згідно з означенням простору $H^{\varphi}(V)$ виконується рівність $v=u\!\upharpoonright\!V$ для деякого $u\in H^{\varphi}(\Omega)$. Оскільки  $(A,B)u\in\mathcal{H}^{\varphi\varrho^{-2q}}(\Omega,\Gamma)$, то на підставі \eqref{implication} маємо включення $u\in C^{p}(\Omega_{0}\cup\Gamma_{0})$. Звідси $v\in C^{p}(\overline{V})$. Таким чином, $H^{\varphi}(V)\subset C^{p}(\overline{V})$, що тягне за собою умову \eqref{1f15} на підставі \eqref{v}. Зауваження~2 обґрунтоване.

\textbf{\textit{Доведення теореми}~6.} Включення $u\in C^{2q}(\Omega)$ є наслідком умов \eqref{1f16} і \eqref{f18} на підставі теореми 5, у якій покладаємо $p:=2q$, $\varphi:=\varphi_1$, $\Omega_{0}:=\Omega$ і $\Gamma_{0}:=\varnothing$. Включення $u\in C^{m}(U_{\sigma}\cup \Gamma)$ є наслідком умов \eqref{1f16}, \eqref{1f17} і \eqref{f19} на підставі тієї ж теореми, у якій беремо $\Omega_{0}:=U_{\sigma}$ і $\Gamma_{0}:=\Gamma$. Таким чином, розв'язок $u$ класичний. Теорема 6 доведена.

\medskip

\end{document}